\documentclass[12pt,twoside]{article}

\usepackage{amsfonts}
\usepackage{amssymb}
\usepackage{amsmath}
\usepackage{graphicx}

\newcommand{\ignore}[1]{}

\def\@begintheorem#1#2{\par\bgroup{\sc #1\ #2. }\it\ignorespaces}
\def\@opargbegintheorem#1#2#3{\par\bgroup{\sc #1\ #2\ (#3). } \it\ignorespaces}
\def\@endtheorem{\egroup}
\newtheorem{theorem}{Theorem}[section]
\newtheorem{corollary}[theorem]{Corollary}
\newtheorem{lemma}[theorem]{Lemma}
\newtheorem{example}[theorem]{Example}
\newtheorem{proposition}[theorem]{Proposition}
\newtheorem{definition}[theorem]{Definition}
\newcommand{\bt}[1]{\begin{theorem}\label{#1}}
\newcommand{\bc}[1]{\begin{corollary}\label{#1}}
\newcommand{\bl}[1]{\begin{lemma}\label{#1}}
\newcommand{\be}[1]{\begin{example}\label{#1}}
\newcommand{\bp}[1]{\begin{proposition}\label{#1}}
\newcommand{\ba}[1]{\begin{algorithm}\rm\label{#1}}
\newcommand{\bd}[1]{\begin{definition}\rm\label{#1}}{\normalsize }
\newcommand{\bpr}{\noindent {\em Proof. }}
\newcommand{\et}{\end{theorem}}
\newcommand{\ec}{\end{corollary}}
\newcommand{\el}{\end{lemma}}
\newcommand{\ee}{\end{example}}
\newcommand{\ep}{\end{proposition}}
\newcommand{\ed}{\end{definition}}
\newcommand{\epr}{{\ \vbox{\hrule\hbox{%
\vrule height1.3ex\hskip0.8ex\vrule}\hrule}}\\\par}

\def\R{\mathbb{R}}
\def\Z{\mathbb{Z}}

\def \supp {{\rm supp}}
\def \lexmin {{\rm lexmin}}

\def\i{\wedge}
\def\u{\vee}

\begin{document}

\title{\bf The Unimodular Intersection Problem}

\author{
Volker Kaibel
\thanks{\small Otto-von-Guericke Universit\"at, Magdeburg, Germany. Email: kaibel@ovgu.de}
\and
Shmuel Onn
\thanks{\small Technion - Israel Institute of Technology, Haifa, Israel.
Email: onn@ie.technion.ac.il}
\and
Pauline Sarrabezolles
\thanks{\small Ecole des Ponts, Paris, and
Technion, Haifa. Email: pauline.sarrabezolles@gmail.com}
}

\date{}

\maketitle

\begin{abstract}
We show that finding minimally intersecting $n$ paths from
$s$ to $t$ in a directed graph or $n$ perfect matchings
in a bipartite graph can be done in polynomial time.
This holds more generally for unimodular set systems.


\end{abstract}

\section{Introduction}

Let $G$ be a directed graph on $d$ edges with two distinct vertices $s$
and $t$ and let $n$ be a positive integer. Given $n$ paths from $s$ to $t$,
an edge is {\em critical} if it is used by all paths. We wish to
find $n$ paths from $s$ to $t$ with minimum number of critical edges.
One motivation for this problem, adapted from \cite{AENYZ,OSZ}, is as follows.
We need to make a sensitive shipment from $s$ to $t$. In the planning stage,
$n$ paths are chosen and prepared. Then, just prior to the actual shipment,
one of these paths is randomly chosen and used. An adversary, trying to harm
the shipment and aware of the prepared paths but not of the actual path
that is finally chosen, will try to harm a critical edge, used by all paths.
Therefore, we protect each critical edge with high cost, and so our goal
is to choose paths with minimum number of critical edges.

Here we study and solve a more difficult problem, described as follows.
Following \cite{AENYZ}, for $1\leq k\leq n$, call an edge {\em $k$-vulnerable}
if it used by at least $k$ of the $n$ paths. In \cite{OSZ} it was shown that
finding $n$ paths with minimum number of $2$-vulnerable edges is NP-hard, and in
\cite{AENYZ} it was shown that for fixed $n$ and $k$, finding $n$ paths
with minimum number of $k$-vulnerable edges can be done in polynomial time.

Here, we assume that $n$ is variable, and want to find $n$ paths with
{\em lexicographically minimal intersection}, that is, which first of all
minimize the number of $n$-vulnerable edges, then of $(n-1)$-vulnerable edges,
and so on. More precisely, given $n$ paths, define their {\em vulnerability vector}
to be $f=(f_1,\dots,f_n)$ with $f_k$ the number of $k$-vulnerable edges.
Then $n$ paths with vector $f$ are better then $n$ paths with vector $g$
if the last nonzero entry of $g-f$ is positive, denoted $f\prec g$. In particular, $n$
paths with lexicographically minimal vector have a minimum number of critical edges.

Let $A$ be the vertex-edge incidence matrix of $G$ and let $b$ be the vector in the
vertex space with $b_s=-1$, $b_t=1$, and $b_v=0$ for any other vertex $v$.
For $x\in\R^d$ let $|x|:=\sum_{i=1}^d |x_i|$.
For $x^1,\dots,x^n\in\R^d$ let $x^1\i\cdots\i x^n\in\R^d$ and $x^1\u\cdots\u x^n\in\R^d$
be the entry-wise minimum and maximum of $x^1,\dots,x^n$,
and define $f(x^1,\dots,x^n)$ by
$$f_k(x^1,\dots,x^n)\ :=\ |\u\{x^{i_1}\i\cdots\i x^{i_k}\ :\
1\leq i_1<\cdots<i_k\leq n\}|\ ,\quad k=1,\dots, n\ .$$
So for $x^1,\dots,x^n\in\{0,1\}^d$ we have that
$f_k(x^1,\dots,x^n)$ is the number of \emph{$k$-vulnerable components} $1\leq i\leq d$, 
in which at least $k$ of the vectors $x^1,\dots,x^n$ have a one-entry.
Then our problem is the following nonlinear integer program defined by $A,b,n$:
\begin{equation}\label{intersection_problem}
\lexmin\{f(x^1,\dots,x^n)\ :\ x^k\in\{0,1\}^d\,,\ Ax^k=b\,,\ k=1,\dots,n\}\ .
\end{equation}
In this article we solve this problem, even in the more general context where $A$ is an arbitrary totally unimodular matrix. We establish the following theorem.
\bt{Main}
The lexicographic problem (\ref{intersection_problem})
over any totally unimodular matrix $A$, any integer right-hand side $b$, and any positive integer $n$,
is polynomial time solvable.
\et

Another consequence of this theorem is the following. Let $G$ be a bipartite graph
on $d$ edges and let $n$ be a positive integer. Given $n$ perfect matchings in $G$,
an edge is $k$-vulnerable if it is used by at least $k$ matchings. The problem of finding
$n$ perfect matchings with lexicographically minimal intersection is again a special case
of our theorem with $A$ the vertex-edge incidence matrix of $G$ and $b$ the vector in the
vertex space with $b_v=1$ for every vertex $v$. So it can also be solved in polynomial time.

\vskip.2cm
In the next section we prove our theorem, proving on the way a core lemma which is of interest
in its right. In the last section we discuss possibilities of some further generalizations, 
solve the simpler problem of minimizing $f_n(x^1,\dots,x^n)=|x^1\i\cdots\i x^n|$ 
only more efficiently, and describe some NP-complete variants.

\section{Proof}
\label{proof}

We begin with an algorithmic version of a decomposition result of \cite{BT}.

\begin{lemma}\label{lemma}
There is a polynomial time algorithm that given any positive integer $n$, 
any totally unimodular matrix $A$, any integer vector $b$, and any vector\break 
$x\in\{0,1,\dots,n\}^d$ satisfying $Ax=nb$, finds vectors $x^1,\dots,x^n\in\{0,1\}^d$ satisfying
$$\sum_{k=1}^n x^k=x\,,\quad Ax^k=b\,,\ k=1,\dots,n\ .$$
\end{lemma}

\vskip.2cm
\bpr
We use induction on $n$. For $n=1$ take $x^1=x$.
Next consider $n>1$ and the system
$0\leq\lfloor{1\over n} x_j\rfloor\leq y_j\leq
\lceil{1\over n} x_j\rceil\leq 1$, $j=1,\dots,d$, $Ay=b$.
Since ${1\over n} x$ is a real solution and $A$ is
totally unimodular, there is also an integer solution $x^n\in\{0,1\}^d$, which can be
found in polynomial time by linear programming. Set $z:=x-x^n$. Then $Az=(n-1)b$,
and for all $j$, if $x_j=0$ then $x^n_j=0$ and if $x_j=n$ then $x^n_j=1$, so
$z\in\{0,1,\dots,n-1\}^d$. Therefore, by induction, we can find in polynomial time
vectors $x^1,\dots,x^{n-1}\in\{0,1\}^d$ satisfying $\sum_{k=1}^{n-1} x^k=z$ and
 $Ax^k=b$ for all $k$. Then $\sum_{k=1}^n x^k=z+x^n=x$
and the induction follows, completing the proof.
\epr
Next, we make the following definition, which is central to what follows. 
The \emph{compression} of $x=(x^1,\dots,x^n)\in\{0,1\}^{dn}$ 
is $\bar{x}=(\bar{x}^1,\dots,\bar{x}^n)\in\{0,1\}^{dn}$ defined by
$$\bar{x}^k\ =\ \u\{x^{i_1}\i\cdots\i x^{i_k}\ :\
1\leq i_1<\cdots<i_k\leq n\}\ ,\quad k=1,\dots,n\ .$$
Thus, the vulnerability vector of $x=(x^1,\dots,x^n)\in\{0,1\}^{dn}$
is $f(x)=(|{\bar x}^1|,\dots,|{\bar x}^n|)$.

Call $c=(c^1,\dots,c^n)\in\R^{dn}$ {\em nondecreasing} if $c^1\leq\cdots\leq c^n$.
For $c,x\in\R^{dn}$ put $cx:=\sum_{k=1}^nc^kx^k$. 
We have the following lemma which is of interest in its own right.   

\begin{lemma}\label{core}
There is a polynomial time algorithm that given any positive integer $n$,
any totally unimodular matrix $A$, any integer vector $b$, and any nondecreasing 
$c=(c^1,\dots,c^n)\in\Z^{dn}$, solves the following nonlinear integer programming problem,
\begin{equation}\label{problem}
    \min\{c{\bar x}\ :\  Ax^k=b\,,\ k=1,\dots,n\,,\ x=(x^1,\dots,x^n)\in\{0,1\}^{dn}\}\ .
\end{equation}
\end{lemma}

Note that the objective $c{\bar x}$ is nonlinear since $c$ acts on the compression of $x$. 

\vskip.2cm\bpr
First, we need a couple of properties of the compression. Note that it satisfies 
$$\bar{x}^k_i\ =\
            \begin{cases}
                1 & \text{if } k\leq\sum_{l=1}^n x^l_i \\
                0 & \text{otherwise}
            \end{cases}\,,\quad k=1,\dots,n\,,\quad i=1,\dots,d\ .$$
Therefore, for all $x,y,z\in\{0,1\}^{dn}$, we have
\begin{equation}\label{compression1}
\sum_{k=1}^n {\bar x^k}\ =\ \sum_{k=1}^n x^k\ ,
\end{equation}
\begin{equation}\label{compression2}
\mbox{and}\ ,\quad\mbox{if}\quad\sum_{k=1}^n y^k\ =\ \sum_{k=1}^n z^k\ ,
\quad\mbox{then}\quad {\bar y}\ = \ {\bar z}\ .
\end{equation}

Now, consider the linear integer program
\begin{equation}\label{relaxation}
    \min\{cx\ :\  A\sum_{k=1}^n x^k=nb\,,\ x=(x^1,\dots,x^n)\in\{0,1\}^{dn}\}\ .
\end{equation}
As $[A,\dots,A]$ is totally unimodular since $A$ is, and $b$ is integer, we can solve 
program \eqref{relaxation} in polynomial time using linear programming. 
If it is infeasible then so is program (\ref{problem}) and we are done. 
So assume we obtain an optimal solution $y\in\{0,1\}^{dn}$ to \eqref{relaxation}. 
By Lemma~\ref{lemma} applied to $\sum_{k=1}^n y^k$ we can find in polynomial time $z\in\{0,1\}^{dn}$ with
$$\sum_{k=1}^n z^k\ =\ \sum_{k=1}^n y^k \,,\quad Az^k=b\,,\ k=1,\dots,n\ .$$
Clearly $z$ is feasible in (\ref{problem}). We claim it is also optimal in \eqref{problem}. 
To prove this, we consider any $x$ which is feasible in (\ref{problem}), 
and prove that the following inequality holds,
$$c{\bar z}\ =\ c{\bar y}\ \leq\ cy\ \leq\ c{\bar x}\ .$$
Indeed, the first equality follows since $\sum_{k=1}^n z^k\ =\ \sum_{k=1}^n y^k$ and therefore
${\bar z}={\bar y}$ by \eqref{compression2}. The middle inequality follows since $c$ is 
nondecreasing. The last inequality follows since $\sum_{k=1}^n {\bar x^k}\ =\ \sum_{k=1}^n x^k$ 
by \eqref{compression1} and therefore ${\bar x}$ is feasible in \eqref{relaxation}.
\epr

\vskip.2cm
We now proceed with the proof of Theorem \ref{Main}. Define a nondecreasing $c\in\Z^{dn}$,
$$c^k_i\ :=\ (d+1)^{k-1}\,,\quad k=1,\dots,n\,,\quad i=1,\dots,d\ .$$
Consider any $x,y\in\{0,1\}^{dn}$, with vulnerability vectors
$f=(|{\bar x}^1|,\dots,|{\bar x}^n|)$ and $g=(|{\bar y}^1|,\dots,|{\bar y}^n|)$.
Suppose $f$ is lexicographically smaller than $g$, that is, $f\prec g$. Let $r$ be the largest 
index such that $|{\bar x}^r|\neq |{\bar y}^r|$. Then $|{\bar y}^r|\geq |{\bar x}^r|+1$. We then have
\begin{eqnarray*}
c{\bar y}-c{\bar x} 
&   =  & \sum_{k=1}^n(d+1)^{k-1}\left(|{\bar y}^k|-|{\bar x}^k|\right) \\
& \geq & \sum_{k<r}(d+1)^{k-1}\left(|{\bar y}^k|-|{\bar x}^k|\right)+(d+1)^{r-1} \\
& \geq & (d+1)^{r-1} -\sum_{k<r}d(d+1)^{k-1}\ >\ 0\ .
\end{eqnarray*}
Thus, an optimal solution $x^1,\dots,x^n$ for \eqref{problem}
is also optimal for our original problem (\ref{intersection_problem}),
and can be found by Lemma \ref{core} in polynomial time . Theorem \ref{Main} follows.
\epr

\section{Remarks}

Gijswijt~\cite{Gij05} established a decomposition result like the one we recalled in
Lemma~\ref{lemma} for the more general case of $A$ being only \emph{nearly totally unimodular},
i.e., $A=A'+ca^T$ with the matrix arising from adding row $a^T$ to $A'$ being totally unimodular
and the vector $c$ being integer. Using this result, one can, for instance, see that by the algorithm
described in Section~\ref{proof} the lexicographically minimal intersection problem cannot only be
solved in polynomial time for structures like $s$-$t$ paths or matchings in bipartite graphs,
but also, for example, for stable sets in circular arc graphs.

\vskip.2cm
Next we discuss the simpler problem of minimizing
$f_n(x^1,\dots,x^n)=|x^1\i\cdots\i x^n|$ only, which in the paths problem
is the number of critical edges, that is, the problem
\begin{equation}\label{simple_problem}
\min\{|x^1\i\cdots\i x^n|\ :\ x^k\in\{0,1\}^d\,,\ Ax^k=b\,,\ k=1,\dots,n\}\ .
\end{equation}
It can be solved in polynomial time as a corollary of Theorem \ref{Main}, but we give here 
a simpler direct proof which encompasses a more efficient algorithm, where we need to solve 
a program in $2d$ variables instead of $nd$ variables. So we have the following.

\bc{Simple}
The intersection problem (\ref{simple_problem})
over any totally unimodular matrix $A$, any integer right-hand side $b$, and any positive integer $n$,
is polynomial time solvable.
\ec
\bpr Consider the program
\begin{equation}\label{auxiliary_problem_2}
\min\{|z|\ :\ y\in\{0,1,\dots,n-1\}^d\,,\ z\in\{0,1\}^d\,,\ A(y+z)=nb\}\ .
\end{equation}
We solve it in polynomial time using linear programming. If it is infeasible then so is 
program (\ref{simple_problem}) and we are done. So assume we obtain an optimal solution $y,z$.

Applying Lemma \ref{lemma} to $x:=y+z$, we can find in polynomial time $x^1,\dots,x^n$ feasible in (\ref{simple_problem})
with $\sum_{k=1}^n x^k=x=y+z$. For $1\leq j\leq d$, if $x^1_j=\cdots= x^n_j=1$ then
$y_j+z_j=\sum_{k=1}^n x^k_j=n$ and $y_j\leq n-1$ so $z_j=1$. Therefore $x^1\i\cdots\i x^n\leq z$.
We claim the $x^k$ are optimal for (\ref{simple_problem}).
Suppose indirectly ${\tilde x}^1,\dots,{\tilde x}^n$ form a better solution.
Then ${\tilde z}:={\tilde x}^1\i\cdots\i {\tilde x}^n$ and ${\tilde y}:=\sum_{k=1}^n{\tilde x}^k-{\tilde z}$
are feasible in (\ref{auxiliary_problem_2}). But then
$$|{\tilde z}|\ =\ |{\tilde x}^1\i\cdots\i {\tilde x}^n|\ <\ |x^1\i\cdots\i x^n|\ \leq\ |z|\ ,$$
which is a contradiction. So $x^1,\dots,x^n$ are optimal and the corollary follows.
\epr
Finally, we point out two variants of problem (\ref{simple_problem}) where 
even deciding if the optimal value in the problem is $0$ is NP-complete,
even for fixed $n=2$. First, if for each $k$ we have a system $Ax^k=b^k$ with the
same matrix $A$ but its own right-hand side $b^k$, then, this is hard even if $A$
is the vertex-edge incidence matrix of a directed graph $G$. Indeed, for $k=1,2$,
let $s_k$ and $t_k$ be distinct vertices in $G$ and let $b^k$ be the vector forcing
$x^k\in\{0,1\}^d$ satisfying $Ax^k=b^k$ to be a path from $s_k$ to $t_k$.
Then the optimal objective function value of problem (\ref{simple_problem})
is $0$ if and only if for $k=1,2$ there are paths from $s_k$ to $t_k$ which are
edge-disjoint, which is a classical NP-complete problem to decide. Second, if $A$ is
not totally unimodular, then this is hard even if $A$ has a single row.
Indeed, let $a_1,\dots,a_d$ be given positive integers with even sum, and consider
the NP-complete problem of deciding if there is a partition, that is,
$I\subseteq\{1,\dots,d\}$ with $\sum_{i\in I}a_i=\sum_{i\notin I}a_i$.
Let $A:=[a_1,\dots,a_d]$, let $b:={1\over2}\sum_{i=1}^da_i$, and let $n:=2$.
Then, if $I$ provides a partition, then $x^1,x^2$ defined by $x^1_i:=1-x^2_i:=1$ if $i\in I$
and $x^1_i:=1-x^2_i:=0$ if $i\notin I$, are feasible in (\ref{intersection_problem})
with value $0$, and if $x^1,x^2$ are feasible in (\ref{intersection_problem}) with value $0$,
then $I:=\supp(x^1):=\{i:x^1_i\neq 0\}$ provides a partition.

\section*{Acknowledgments}

The research of the first author was partially supported by Deutsche Forschungsgemeinschaft (KA 1616/4-1).
The research of the second author was partially supported by a VPR Grant
at the Technion and by the Fund for the Promotion of Research at the Technion.
The research of the third author was partially supported by a
Bourse d'Aide \`a la Mobilit\'e Internationale from Universit\'e Paris Est.

\end{document}